\documentclass[12pt, a4paper]{amsart}
\usepackage[T1]{fontenc}
\usepackage{amsfonts}
\usepackage{listings}
\usepackage{color}
\usepackage{fancyvrb}
\usepackage{amsmath}
\usepackage{mathrsfs}
\usepackage{amsthm}
\usepackage{graphicx}
\usepackage{enumerate}
\usepackage{bbm}
\usepackage{euscript}
\usepackage{yhmath}
\usepackage{MnSymbol}
\usepackage{parskip}
\usepackage[colorlinks=true,
            linkcolor=blue,  
            citecolor=blue,  
            urlcolor=blue]{hyperref}
\usepackage[a4paper, margin=1in]{geometry}
\newcommand{\R}{\mathbb{R}}

\newcommand{\ep}{\varepsilon}

\newcommand{\nab}{\nabla}

\newcommand{\p}{\partial}
\newcommand{\LL}{\mathrm{L}}
\newcommand{\WW}{\mathrm{W}}
\newcommand{\CC}{\mathrm{C}}
\newcommand{\GG}{\mathrm{G}}
\newcommand{\YY}{\mathrm{Y}}

\newcommand{\NN}{\mathrm{N}}

\newcommand{\dd}{\mathrm{d}}
\newcommand{\BMO}{\mathrm{BMO}}

\newtheorem{theorem}{Theorem}[section]
\newtheorem{lemma}[theorem]{Lemma}
\newtheorem{prop}[theorem]{Proposition}

\theoremstyle{definition}

\DeclareMathOperator{\Supp}{supp}
\DeclareMathOperator{\dist}{dist}
\DeclareMathOperator{\rank}{rank}

\DeclareMathOperator{\Div}{div}

\title{Regularity of minimisers for the $p$-distance functional}
\author{Aidan Strong}
\begin{document}
\baselineskip=18pt

\maketitle
\begin{abstract}
We consider regularity of minimisers of the $p$-distance functional $\mathcal{F}[u]=\int_\Omega \mathrm{dist}(Du,K)^p \dd x$ where $K \subset \R^{N \times n}$ is a compact set. We distinguish the scalar and vectorial cases, and consider geometric conditions on $K$ where regularity can be improved.
\end{abstract}

\setcounter{tocdepth}{1}
\begingroup
\small
\tableofcontents
\endgroup

\section{Introduction}
In this paper we study the regularity of minimisers of \begin{equation}\label{eqn distance functional}
\mathcal{F}[u]=\int_\Omega \dist(Du,K)^p \dd x
\end{equation}
where $n \geq 2$, $N \geq 1$, $\Omega \subset \R^n$ is a bounded domain, $K \subset \R^{N \times n}$ is a compact set, $1<p<\infty$, and $u \in \WW^{1,p}(\Omega;\R^{N})$. We refer to $\mathcal{F}$ as the $p$-distance functional, and we write $F_K(z)=\dist(z,K)$. There are several motivating reasons to study the regularity of minimisers of this functional. Firstly, in the special case $K=\left\{ 0 \right\}$ our functional is simply the $p$-Dirichlet integral. In this case the regularity of minimisers (which are precisely $p$-harmonic maps) is a classical topic going back to the work of Uhlenbeck \cite{Uhlenbeck1977}, who showed that $u \in \CC^{1,\alpha}_\text{loc}(\Omega;\R^{N})$. From a purely theoretical perspective, the $p$-distance functional is therefore a natural extension of the $p$-Dirichlet integral. Secondly, our functional is a typical example where the integrand vanishes if $Du$ lies in some prescribed set $K$. Functionals with this property naturally arise in non-linear elasticity, where the compact set $K$ can be regarded as the set of preferred deformations. A typical example is the multi-well model $K= \bigcup_{i=1}^mSO(n) U_i$ where $U_i$ are positive-definite symmetric matrices. Related to this, a third motivation is the relation to rigidity problems. Indeed, if the differential inclusion $Du \in K$ holds almost everywhere, then $u$ is obviously a minimiser of \eqref{eqn distance functional}. We can think of minimisers of the distance functional as the maps closest to satisfying the differential inclusion $Du \in K$ over some Dirichlet class. The regularity and rigidity of solutions to differential inclusions are naturally related to the geometry of $K$, and as we will see this is also the case for minimisers of the $p$-distance functional. For example if $K$ is convex, the lack of rank-one connections in $K$ results in regularity for both minimisers and solutions of the differential inclusion. If $K$ is not convex it turns out the distance functional is not weakly lower semicontinuous on $\WW^{1,p}$, so we instead shall consider the relaxation in this case:
\begin{equation}\label{eqn distance functional quasiconvex envelope}
    \mathcal{F}^{qc}[u]=\int_\Omega (F_K^p)^{qc}(Du) \dd x.
\end{equation}
When $K$ is convex, the relaxed functional \eqref{eqn distance functional quasiconvex envelope} obviously agrees with \eqref{eqn distance functional}.
A systematic approach to the regularity of minimisers of \eqref{eqn distance functional} and \eqref{eqn distance functional quasiconvex envelope} was initiated in \cite{DolzmannKristensenZhang2013}, where several results and counterexamples were shown in the case $p=2$. In essence, this paper is the continuation of \cite{DolzmannKristensenZhang2013} for the full range $1<p<\infty$.
Let us now state our main results. Firstly, we address the non-convex case.

\begin{theorem}\label{Theorem BMO regularity for distance functional}
Let $1<p<\infty$ and $K \subset \R^{N \times n}$ compact. Suppose that $u \in \WW^{1,p}(\Omega;\R^N)$ is a local minimiser of the relaxed functional \eqref{eqn distance functional quasiconvex envelope}. Then $$Du \in \BMO_\text{loc}(\Omega;\R^{N \times n}).$$
\end{theorem}
The case $p=2$ was already obtained in \cite{DolzmannKristensenZhang2013}. We will give a short proof here relying on more general results from \cite{Strong2026} where $\BMO$ regularity is established for functionals that are suitably asymptotic to $|z|^p$ in a first order sense. As shown in \cite{DolzmannKristensenZhang2013}, one can construct non--Lipschitz minimisers in the setting of Theorem \ref{Theorem BMO regularity for distance functional} in several examples. In particular, there are non-Lipschitz minimisers if $N \geq 2$ when $K=SO(2)$ or $K=\left\{A,-A\right\}$ where $A^TA=I$ and $\rank A >1$. A striking feature of the distance functional compared with $p$-harmonic maps, is that the vectorial and scalar cases behave differently. In the scalar case we have the following.

\begin{theorem}\label{Theorem scalar case regularity non convex distance functional}
Let $1<p<\infty$ and $K \subset \R^n$ compact. Suppose that $u \in \WW^{1,p}(\Omega)$ is a local minimiser of \eqref{eqn distance functional quasiconvex envelope} where $N=1$. Then $u \in \WW^{1,\infty}_\text{loc}(\Omega)$ and
$$\sup_{B_{R/2}}|Du|^p \ \leq C \left(1+\strokedint_{B_R}|Du|^p \dd x \right)$$
for all $B_R \Subset \Omega.$
\end{theorem}
Our proof of Theorem \ref{Theorem scalar case regularity non convex distance functional} is in the spirit of \cite{Fonseca1997, Fonseca2002}, where one uses Moser iteration combined with suitable approximations to deal with non-differentiability and degeneracy. We note also that in the scalar case, the quasiconvex envelope is just the convex envelope which in turn is given by $F_{\mathrm{co}K}^p(z)$. This means that for the scalar case, we may assume $K$ is convex. The results of  \cite{Fonseca1997, Fonseca2002} allow one to prove Lipschitz regularity for a large class of degenerate integrands. In particular, Lipschitz regularity holds whenever $F$ can be approximated by a sequence $F_j$ with uniform $p$-growth bounds, which satisfy both $$\ep_j (1+|z|^2)^\frac{p-2}{2}|\lambda|^2 \leq F_j''(z) \leq \ep_j^{-1}(1+|z|^2)^\frac{p-2}{2},$$ 
$$F_j''(z) \geq \nu(1+|z|^2)^{\frac{p-2}{2}} \text{ for }|z| \geq M,$$
where $\nu$ and $M$ are independent of $j$. For general convex sets $K$, this second condition need not hold for any approximating sequence. To get around this, we carry out a modified Moser iteration which is adapted to the convexity estimates of $F_K^p$.  It turns out that Theorem \ref{Theorem scalar case regularity non convex distance functional} is sharp. In the scalar case, if $K$ contains two distinct points, by considering a laminate one can always construct a minimiser which is not $C^1$. On the other hand, if $K$ is a singleton then we know from the theory of $p$-harmonic maps that $u \in C^{1,\alpha}_\text{loc}(\Omega)$. More generally, if $K \subset \R^{N \times n}$ then a non $C^1$-minimiser can be constructed whenever $K$ has a rank-one connection. The question of whether Theorem \ref{Theorem scalar case regularity non convex distance functional} holds in the vectorial setting under the assumption that $K$ is convex was asked in \cite{DolzmannKristensenZhang2013}. We address this in the negative by means of a counterexample.
\begin{theorem}\label{theorem counterexample lipschitz vectorial first one}
There exists a compact convex set $K \subset \R^{9 \times 3}$ and a map $u \in \WW^{1,2}(B_1(0);\R^{9})$, such that $u$ is a local minimiser of $\int_{B_1(0)}F_K^2(Dv) \dd x$, but $u$ is not Lipschitz.
\end{theorem}
Our construction relies on the observation that $u \in \WW^{1,2}$ is a local minimiser of the square distance functional if and only if:
\begin{enumerate}[(i)]
    \item $Du=P+\nu$ where $\nu$ is divergence free and $P \in K$.
    \item $\nu(x) \cdot (k-P(x)) \leq 0$ for all $k \in K$ and $x \in \R^n$.
\end{enumerate}
Given a map $u$, writing it as $Du=P+\nu$ where $\nu$ is divergence free, we can take $K$ as the convex hull of the image of $P$, so that the supporting inequality $(ii)$ becomes $\nu(x) \cdot (P(y)-P(x)) \leq 0$ for all $x,y \in \R^n$. Therefore, it suffices to find some non-Lipschitz map $u \in \WW^{1,2}$ where this supporting inequality holds for some divergence free splitting of $Du$. The dimensions $(n,N)=(3,9)$ are not optimal, but allow us to give a fairly straightforward construction. With the assistance of AI tools, we have been able to produce a counterexample in the case $(n,N)=(2,3)$, leaving only the case $(n,N)=(2,2)$ open. However, we intend to report on this elsewhere. The final problem we study is whether one can obtain better regularity if $K$ satisfies certain geometric conditions. A compact convex set $K \subset \R^m$ satisfies SBC\footnote{SBC stands for "supporting ball condition"} if for every $z \in \partial K$ there exists $R_z >0$ and a closed ball $\overline B_{R_z}(c_z)$ such that $z \in \partial \overline B_{R_z}(c_z)$ and $K \subset \overline B_{R_z}(c_z).$ $K$ satisfies uniform SBC if $R_z$ can be chosen to be independent of $z$. The following is a generalisation of a result from  \cite{DolzmannKristensenZhang2013} which dealt with $p=2$. 

\begin{theorem}\label{theorem lipschitz regularity for distance functional under supporting ball condition}
Let $K \subset \R^{N \times n}$ be a compact convex set satisfying uniform SBC, and let $1<p<\infty$. Suppose that $u \in \WW^{1,p}(\Omega;\R^N)$ is a local minimiser of \eqref{eqn distance functional}. Then $u \in \WW^{1,\infty}_\text{loc}(\Omega;\R^N)$ and $$\sup_{B_{R/2}}|Du|^p \ \leq C \left(1+\strokedint_{B_R}|Du|^p \dd x \right)$$
for all $B_R \Subset \Omega.$
\end{theorem}
The proof is very different to the scalar case, and we once more use some general results from \cite{Strong2026} concerning functionals which are suitably asymptotic to $|z|^p$, but this time in a second order sense. It is an interesting problem to characterise all compact convex sets $K \subset \R^{N \times n}$ for which Lipschitz regularity holds. For example, if $p=2$ then Proposition \ref{prop vectorial regularity decouples} gives Lipschitz regularity when $K=[-1,1]^{N \times n}$, which clearly does not satisfy uniform SBC. This suggests that it is not simply flatness of $K$ that is an obstruction to Lipschitz regularity. To conclude, we consider some special examples where further regularity holds in the case $p=2$. Firstly, we consider the case when $K$ has no rank-one connections. Via a difference quotient argument, this leads to $W^{2,2+\ep}_\mathrm{loc}$ regularity of minimisers. Finally, we consider the case $K=[-A,A]$ is a line-segment. In the case $\rank A \geq 2$, this reduces to a scalar problem where De-Giorgi--Nash--Moser theory can be applied to conclude $Du \in C^{0,\alpha}_\mathrm{loc}$.
\section{Preliminaries}
\subsection{Notation}
Throughout this paper, $\Omega \subset \R^n$ is an open bounded set where $n \geq 2$. We consider maps $u \in \WW^{1,p}(\Omega;\R^{N})$ where $1<p<\infty$ and write $\R^{N \times n}$ for the space of $N\times n$ matrices. We equip $\mathbb R^{N\times n}$ with the standard
Euclidean inner product denoted $z \cdot w$, and denote the associated norm by $|z|$. We write
$p'=p/(p-1)$ for the conjugate exponent of $p$. For $x_0\in\mathbb R^n$ and $R>0$, we denote by
$B_R(x_0)$ the open ball of radius $R$ centred at
$x_0$, and as is standard we simply write $B_R$ when it is convenient. For an integrable function $f$, we denote the mean over the ball $B_R(x_0)$ by $$(f)_R=(f)_{B_R(x_0)}=\strokedint_{B_R}f \dd x = \frac{1}{|B_R|}\int_{B_R}f \dd x.$$
As is standard, $C,c>0$ may refer to constants changing from line to line. We call a constant $C$ universal if it only depends on the data $n,N,p$. We write $A\lesssim B$ if $A\leq CB$ for some constant $C>0$, and $A \approx  B$ if $A\lesssim B$ and $B\lesssim A$. Given an integrand $F:\R^{N\times n}\rightarrow  \R$ with $p$-growth, we say that $u\in \WW^{1,p}_{\mathrm{loc}}
(\Omega;\mathbb R^N)$ is a local minimiser of
$\int_\Omega F(Dv)\dd x$ if
$$\int_{\Omega'} F(Du) \dd x
\leq \int_{\Omega'} F(Du+D\varphi)\dd x$$
for every $\Omega' \Subset \Omega$ and
$\varphi\in W^{1,p}_0(\Omega';\R^N)$.
\subsection{Auxiliary functions}
Define for $z \in \R^{N \times n}$
$$A_p(z)=|z|^{p-2}z,$$
$$V_p(z)=|z|^{\frac{p-2}{2}}z.$$
Note that $|z|^p$ is a $C^1$ function whose derivative is $pA_p(z).$ We retrieve the following from \cite{Strong2026}.
\begin{lemma}\label{lemma equivalences of auxiliary V functinos}
Let $1<p<\infty$ and $z,w \in \R^{N \times n}$. The following quantities are comparable up to universal constants:
\begin{enumerate}[(i)]
    \item $|V_p(z)-V_p(w)|^2$
    \item $(|z|+|w|)^{p-2}|z-w|^2$
    \item $(A_p(z)-A_p(w)) \cdot (z-w)$
    \item $|A_p(z)-A_p(w)| |z-w|$
    \item $ |z|^p -|w|^p -pA_p(w) \cdot (z-w)$
    \item $|z-w|^2 \int_0^1|tz+(1-t)w|^{p-2} \dd t$.
\end{enumerate}
\end{lemma}
The following technical lemma is useful to understand the asymptotic nature of the $p$-distance function.

\begin{lemma}\label{lem: Bregman divergence for V function}
Let $1<p<\infty$, and suppose that $G \in C(\R^m;\R^l)$ has a locally  Lipschitz derivative outside a sufficiently large ball, and satisfies
$$|G(z)| \leq L(1+|z|^{p-1})$$ for all $z$, and 
$$|G'(z)| \lesssim |z|^{p-2}, \qquad |G''(z)| \lesssim |z|^{p-3}$$
for almost all sufficiently large $|z|$. Then there exist $C,M>0$ such that
$$|G(z)-G(a)-G'(a) (z-a)| \leq \frac{C}{|a|}|V_p(z)-V_p(a)|^2$$
for all $z \in \R^m$ and for all $|a| \geq M$. In particular, for $|w| \geq 1$ we have
$$|A_p(z)-A_p(w)-A_p'(w)(z-w)| \lesssim \frac{1}{|w|}(|z|^p-|w|^p-pA_p(w)\cdot (z-w)).$$
\end{lemma}

\begin{proof}
We split the inequality into two cases.\\
\underline{Case 1:} $|z-a| \leq \frac{|a|}{2}$. 
The triangle inequality gives $\frac{1}{2}|a| \leq |z| \leq \frac{3}{2}|a|$, and in particular $|a+t(z-a)| \approx |a|$ for $0<t<1$. So as long as $|a|$ is sufficiently large, we may apply the Fundamental Theorem of Calculus and the bound on $G''$ to deduce
\begin{align*}
|G(z)-G(a)-G'(a)(z-a)| &=\left|\int_0^1G'(a+t(z-a))-G'(a)\dd t(z-a) \right|\\
& \leq |a|^{p-3}|z-a|^2 \\
& \approx \frac{(|z|+|a|)^{p-2}|z-a|^2}{|a|}\\
& \approx \frac{|V_p(z)-V_p(a)|^2}{|a|}.
\end{align*}
\underline{Case 2:} $|z-a| > \frac{|a|}{2}.$
This time the triangle inequality gives $|z-a| \approx |z|+|a|$. We then crudely estimate for large $|a|$
\begin{align*}|G(z)-G(a)-G'(a) (z-a)|&
\leq |G(z)|+|G(a)|+|G'(a)||z-a|\\
&\lesssim 1+|z|^{p-1}+|a|^{p-1}+|a|^{p-2}|z-a|\\
& \lesssim \frac{|a|+(|z|+|a|)^{p}+|a|^{p-1}|z-a|}{|a|}\\
& \approx \frac{(|z|+|a|)^{p}}{|a|}\\
& \approx \frac{(|z|+|a|)^{p-2}|z-a|^2}{|a|}\\
& \approx \frac{|V_p(z)-V_p(a)|^2}{|a|}.
\end{align*}
The final part concerning the $A_p$ function then follows from Lemma \ref{lemma equivalences of auxiliary V functinos}.
\end{proof}

\subsection{Distance function}
Throughout this section, we let $K \subset \R^{N \times n}$ be a non-empty compact set. We define the distance function $$F_K(z)=\dist(z,K).$$ Associated to $K$, we define $$\pi(z)=\left\{k \in K: F_K(z)=|z-k| \right\},$$
$$g(z)=\max_{k \in K}\left( z\cdot k-\frac{1}{2}|k|^2\right).$$
We call $\pi$ the metric projection associated to $K$. The function $g$ is sometimes called the Asplund function, as it was introduced in \cite{Asplund1969} in the study of Chebyshev sets. Since $g$ is a supremum of affine functions, it is always convex. Moreover, from the definition one sees that $$g(z)=\frac{1}{2}(|z|^2-F_K^2(z))$$
and $g$ is globally Lipschitz continuous. By Danskin's theorem (see \cite[Proposition 4.5.1]{Bertsekas2003}), we compute the subdifferential of $g$ as
$$\partial g(z)=\mathrm{conv}(\pi(z))$$
where $\mathrm{conv}$ denotes the convex hull. Since the quantity $z-\pi(z)$ appears in many calculations, it is useful to introduce the notation $$q(z)=z-\pi(z).$$ The convexity of $K$ is intrinsically tied to properties of $F_K$. The following is a restatement of the fact that every Chebyshev set in Euclidean space is convex. 

\begin{prop}\cite{Fletcher2014}
Let $K \subset \R^{N \times n}$ be compact. The following are equivalent.
\begin{enumerate}[(i)]
    \item $K$ is convex.
    \item $F_K$ is convex.
    \item $\pi(z)$ is a singleton for every $z \in \R^{N \times n}$.
    \item $g$ is differentiable at every $z \in \R^{N \times n}$. 
\end{enumerate}
Furthermore, in this case $ g'=\pi$ and $g \in C^{1,1}(\R^{N \times n})$.
Moreover, in this case $$\Vert \pi \Vert_{\text{Lip}} \leq 1$$ and $$q(z) \cdot (k-\pi(z)) \leq 0$$
for all $k \in K$ and $z \in \R^{N \times n}$.
\end{prop}
Assuming still that $K$ is convex, this last inequality is just the statement that $q(z) \in N_K(\pi(z))$ where $$N_K(w) = \left\{ \nu \in \R^{N \times n}: \nu \cdot (k-w) \leq 0 \text{ for all }k\in K \right\}$$
denotes the normal cone at $w \in K$. In fact, it is not hard to show we have
$$\pi^{-1}(w)=N_K(w)+w$$
for $w \in K$, and we note in particular $\pi(\pi(z)+tq(z))=\pi(z)$ for all $t \geq 0$.  Let us collect some monotonicity inequalities, which just follow from the fact that $q(z) \in N_K(\pi(z))$ and Lemma \ref{lemma equivalences of auxiliary V functinos}. 
\begin{lemma}\label{lemma monotonicity estimates for sigma}
Let $K \subset \R^{N \times n}$ be compact and convex, $1<p<\infty$. Write $\sigma(z)=A_p(q(z))$. For $z,w \in \R^{N \times n}$ we have:
\begin{enumerate}[(i)]
    \item  $(\sigma(z)-\sigma(w))\cdot  (\pi(z)-\pi(w)) \geq 0$
    \item $(\sigma(z)-\sigma(w)) \cdot (z-w) \gtrsim |V_{p'}(\sigma(z))-V_{p'}(\sigma(w))|^2$, where $p'=\frac{p}{p-1}.$
\end{enumerate}
\end{lemma}
\begin{proof}
From the above we know $\sigma(z)=A_p(q(z)) \in N_K(\pi(z))$. This gives us both $\sigma(z) \cdot (\pi(w)-\pi(z)) \leq 0$ and $\sigma(w) \cdot (\pi(z)-\pi(w)) \leq 0$. Summing these gives $(i)$. Then $(ii)$ follows by writing $\pi(z)-\pi(w)=z-w -(q(z)-q(w))$, and applying Lemma \ref{lemma equivalences of auxiliary V functinos} together with $V_p=V_{p'} \circ A_p$. 
\end{proof}

The above says that convexity of $K$ is equivalent to convexity of $F_K^p$ for $1 \leq p<\infty$. From the point of view of Calculus of Variations, it is important to understand whether $F_K^p$ can be quasiconvex. It turns out that even rank-one convexity of the distance function is equivalent to convexity. The following can be found in \cite{Zhang1999}, however we give an elementary proof for completeness.

\begin{prop}
Suppose that $K \subset \R^{N \times n}$ is not convex. Then for every $1 \leq p <\infty$, $F_K^p$ is not rank-one convex and hence not quasiconvex. 
\end{prop}

\begin{proof}
Since $K$ is not convex, there exists some
$z_0\in \R^{N\times n}\setminus K$ such that $\pi(z_0)$ is not a singleton. Choose
distinct $w_1,w_2\in \pi(z_0)$ and set $r:=F_K(z_0)=|z_0-w_1|=|z_0-w_2|>0.$
Since $w_1\neq w_2$ and rank-one matrices span $\R^{N \times n}$, there exists a rank-one matrix $H\in \R^{N\times n}$ such that $
H \cdot (w_2-w_1)>0.$ For $t>0$, using $w_2$ and $w_1$ as competitors, respectively, we have
\[
F_K^p(z_0+tH)\leq |z_0+tH-w_2|^p,
\qquad
F_K^p(z_0-tH)\leq |z_0-tH-w_1|^p.
\]
Since $r>0$, the function $v\mapsto |v|^p$ is differentiable at $z_0-w_i$. Hence, as
$t\searrow0$,
\[|z_0+tH-w_2|^p=r^p+pr^{p-2}t\,H\cdot (z_0-w_2)+o(t),\]
and
\[|z_0-tH-w_1|^p=r^p-pr^{p-2}t\,H\cdot (z_0-w_1)+o(t).\]
Therefore
\[\frac{F_K^p(z_0+tH)+F_K^p(z_0-tH)}{2} \leq r^p+\frac{p}{2}r^{p-2}t\,H\cdot (w_1-w_2)+o(t).\]
But $H\cdot (w_1-w_2)<0$, so for all sufficiently small $t>0$,
\[\frac{F_K^p(z_0+tH)+F_K^p(z_0-tH)}{2}<r^p=F_K^p(z_0).\]
It follows $t \mapsto F_K^p(z_0+tH)$ is not convex, and hence $F_K^p$ is not rank-one convex. Since quasiconvexity implies rank-one convexity,
$F_K^p$ is not quasiconvex.
\end{proof}
Next, we collect some facts about the differentiability of $F_K^p$ and its quasiconvex envelope. We begin with the convex case.
\begin{prop}\cite{Fletcher2014}
Let $K \subset \R^{N \times n}$ be compact. Then $F_K$ is differentiable at $z \in \R^{N \times n} \setminus K$ if and only if $\pi(z)$ is a singleton, and in this case
$$F_K'(z) = \frac{q(z)}{|q(z)|}.$$
In particular, if $K \subset \R^{N \times n}$ is convex then $F_K \in C_\text{loc}^{1,1}(\R^{N \times n} \setminus K).$ Moreover, if $K \subset \R^{N \times n}$ is convex and $1<p<\infty$ then
$$(F_K^p)'(z)=pA_p(q(z)).$$ In particular,
$$F_K^p \in \begin{cases}
C^{1,p-1}(\R^m) & \text{if } 1<p \leq 2 \\
C^{1,1}_{\text{loc}}(\R^m) & \text{if } p>2. \end{cases}$$
\end{prop}
Note that Lemma \ref{lemma monotonicity estimates for sigma} provides the quantitative monotonicity estimate of $(F_K^p)'$. 
When $K$ is not convex, we can represent the derivative of $(F_K^p)^{qc}$ in terms of a homogenous gradient Young measure. We refer for instance to \cite{Rindler} for an overview of Young measures and quasiconvex envelopes.

\begin{prop}\cite{BallKirchheimKristensen2000}
Let $K \subset \R^{N \times n}$ be compact. Then $$(F_K^p)^{qc} \in \begin{cases} C^{1,{p-1}}(\R^{N \times n}) & \text{ if } 1<p\leq 2\\
C^{1,1}_{\text{loc}}(\R^{N \times n}) & \text{ if } p>2.\end{cases}$$ Moreover for each $z \in \R^{N \times n}$ there exists a homogeneous gradient Young measure $\nu^z \in \GG \YY^p_\text{hom}(\R^{N\times n})$ with $\overline{\nu^z}=z$ such that $$(F_K^p)^{qc}(z)=\int_{\R^{N\times n}}(F_K^p)(a) \dd \nu^z(a)$$ where $\Supp \nu^z \subset \left\{a: F_K^p(a)=(F_K^p)^{qc}(a) \right\}$. Moreover, $\pi$ is single valued on $\Supp \nu^z$ and $$\nab (F_K^p)^{qc}(z)=\int_{\R^{N \times n}}pA_p(q(a)) \dd \nu^{z}(a).$$
\end{prop}
One can also show that $(F_K^p)^{qc}=(F_{qc(K)}^p)^{qc}$ where $qc(K)$ denotes the quasiconvex hull of $K$. Here, the quasiconvex hull of $K$ is the set of points $z$ where $(F_K^p)^{qc}(z)=0$, which can be equivalently characterised as the set of barycentres of homogenous gradient Young measures with support in $K$, see \cite{Rindler}. We emphasise it is not the case that $(F_K^p)^{qc}$ equals $F_{qc(K)}^p$. In the scalar case $N=1$, we do in fact have $(F_K^p)^{qc}=(F_K^p)^{**}=F_{\mathrm{conv}K}^p$ where $**$ denotes the convex envelope. For details we refer to \cite{Zhang1998, Zhang1999}. Next, we collect some facts about the growth and asymptotic nature of $(F_K^p)^{qc}$. In the language of \cite{Strong2026}, the key point is that $(F_K^p)^{qc}$ resembles $|z|^p$ in a first order sense at infinity. 

\begin{lemma}\label{Lemma growth conditions of distance functional}
Let $1<p<\infty$ and $K \subset \R^{N \times n}$ be compact.
\begin{enumerate}[(i)]
    \item $\mu |z|^p -C \leq (F_K^p)^{qc}(z) \leq L(1+|z|^p)$ for some $\mu,L>0$
    \item $|(F_K^p)^{qc}(z)-|z|^p| \leq L(1+|z|^{p-1})$
    \item $|\nab (F_K^p)^{qc}(z)-pA_p(z)| \leq L(1+|z|)^{p-2}.$
\end{enumerate}
\end{lemma}
\begin{proof}
From the definition $(|z|-M)_+ \leq F_K(z) \leq |z|+M$ where $M=\max_{w \in K}|w|.$ Raising this to the power $p$, and noting that $|z| \mapsto (|z|-M)_+^p$ is convex, it follows that $$(|z|-M)_+^p \leq (F_K^p)^{qc}(z) \leq (|z|+M)^p.$$ From this we immediately get $(i)$ and $(ii)$. The proof of $(iii)$ is more involved. We define $$m(z)=\int_{\R^{N \times n}}q(a) \dd \nu^z(a)=z-\int_{\R^{N \times n}}\pi(a)d\nu^z(a).$$ Since $\pi$ is bounded, we obviously have $|m(z)-z |\leq L$. In particular using Lemma \ref{lemma equivalences of auxiliary V functinos}, $$|A_p(m(z))-A_p(z)| \leq (|m(z)|+|z|)^{p-2}|m(z)-z| \leq L(1+|z|)^{p-2}.$$
Recalling the formula for $\nab (F_K^p)^{qc}$, it suffices to show that $$\left|\int_{\R^{N \times n}}A_p(q(a)) \dd \nu^z(a) -A_p(m(z))\right| \leq L(1+|z|)^{p-2}.$$
To do this, we use Lemma \ref{lem: Bregman divergence for V function} and the definition of $m(z)$, so that for large $z$,
\begin{align*}
&\left|\int_{\R^{N \times n}}A_p(q(a)) - A_p(m(z)) \dd \nu^z(a) \right|\\
&= \left| \int _{\R^{N \times n}}A_p(q(a))-A_p(m(z))-A_p'(m(z))(q(a)-m(z)) \dd \nu^z(a) \right|\\
& \lesssim \frac{1}{|m(z)|} \int_{\R^{N \times n}}|q(a)|^p-|m(z)|^p-pA_p(m(z)) \cdot (q(a)-m(z)) \dd \nu^z(A)\\
& = \frac{1}{|m(z)|}\left((F_K^p)^{qc}(z) -|m(z)|^p\right).\\
\end{align*}
Using $(ii)$ and the definition of $m$, this is clearly bounded by $L(1+|z|)^{p-2}$ so we are done.
\end{proof}
The last thing we consider is the convexity of $F_K^p$ when $K$ is convex. Recalling that $F_K^p \in C^{1,\min(p-1,1)}_\text{loc}(\R^{N \times n}) \cap \WW_\text{loc}^{2,\infty}(\R^{N \times n}\setminus K)$, it makes sense to consider estimates for the second derivative of $F_K^p$ outside of $K$. In the following we write $A \leq B$ for the Loewner order on positive semi-definite matrices. A fact we will use later is that the trace operator is monotone with respect to this order.

\begin{lemma}\label{lemma inequality of q, pi and convexity of distance function}
Let $K \subset \R^m$ be a compact convex set. At points of differentiability $z \in \R^m \setminus K$ of $\pi$, we have:
\begin{enumerate}[(i)]
    \item $0 \leq \pi' \leq I$ where $\pi'$ is symmetric
    \item $0 \leq q' \leq I$ where $q'$ is symmetric
    \item $q'q=q$
    \item $q \otimes q\leq |q|^2 q'$
    \item $\left(\frac{1}{p}F_K^p \right)''=|q|^{p-2}q'+(p-2)|q|^{p-4}q\otimes q.$
    \item $(p-1)|q|^{p-4}q \otimes q \leq \left(\frac{1}{p}F_K^p\right)''(z)\leq \max(p-1,1)|q|^{p-2}I$.
\end{enumerate}
\end{lemma}
\begin{proof}
$(i)$ follows from the fact that $\Vert\pi  \Vert_\text{Lip} \leq 1$, and $\pi' = g''$ is the hessian of a convex function so is symmetric and positive semi-definite. Since $q'=I-\pi'$, we immediately have $(ii)$. To see $(iii)$, recall that $t \mapsto \pi(\pi(z)+tq(z))$ is constant, so differentiating this at $t=1$ gives $\pi'(z)q(z)=0$, from which $(iii)$ follows. For $(iv)$, we estimate using $(ii)$,$(iii)$ and Cauchy--Schwarz
\begin{align*}
(q\otimes q) \xi \cdot \xi = (q \cdot \xi)^2 = (q'q \cdot \xi)^2 \leq |q'q \cdot q||  q'\xi \cdot \xi|=|q|^2 q'\xi \cdot \xi .
\end{align*}
$(v)$ follows from a straightforward calculation using $(iii)$, noting that $(q \otimes q)q'=q \otimes (q'^Tq)=q \otimes (q'q)=q\otimes q$ where we use that $q'$ is symmetric. Finally $(vi)$ follows by combining $(v)$ and $(iv)$.
\end{proof}

 One can obtain better estimates by assuming a stronger notion of convexity. Recall that a compact convex set $K \subset \R^{N \times n}$ satisfies SBC if for every $z \in \partial K$ there exists $R_z >0$ and a closed ball $\overline B_{R_z}(c_z)$ such that $z \in \partial \overline B_{R_z}(c_z)$ and $K \subset \overline B_{R_z}(c_z).$ $K$ satisfies uniform SBC if $R_z$ can be chosen to be independent of $z$. If $K$ has $C^2$ boundary, then uniform SBC is equivalent to SBC, which in turn is equivalent to the positivity of all principal curvatures. One can also characterise uniform SBC as sets of the form $$K=\bigcap_{\alpha}\overline B_R(c_\alpha)$$ for some family of centres $(c_\alpha)$ and some $R>0$. One can also show that uniform SBC is equivalent to $\mathrm{esssup}_{|z|>r}|\pi'(z)|<1$ for some $r>0$, which in turn is equivalent to $\pi' \leq \frac{R}{R+|q|}I$ at points of differentiability \cite{balashov2024strong}. By similar reasoning to before we therefore have the following.

\begin{lemma}\label{lemma estimates SBC}
Let $K \subset \R^m$ be a compact convex set satisfying uniform SBC with radius $R>0$. At points of differentiability $z \in \R^m \setminus K$ of $\pi$, we have:
\begin{enumerate}[(i)]
    \item $0 \leq \pi'\leq \frac{R}{R+|q|}I$
    \item $\frac{|q|}{R+|q|}I \leq q' \leq I$
    \item $\min\left(p-1, \frac{|q|}{R+|q|}\right)|q|^{p-2}I\leq \left(\frac{1}{p}F_K^p\right)''\leq \max(p-1,1)|q|^{p-2}I$.
\end{enumerate}
\end{lemma}

%I should add a proof later maybe since its actually non trivial.

\section{Lipschitz regularity in the scalar case}
In this section we address the scalar case $N=1$, namely Theorem \ref{Theorem scalar case regularity non convex distance functional}. From the previous section we know that the relaxation of \eqref{eqn distance functional} is given by 
\begin{equation}\label{eqn distance functional relaxed scalar}
\int_{\Omega}F_{\mathrm{conv K}}^p(Du) \dd x
\end{equation}
In particular, if $u \in \WW^{1,p}_\text{loc}(\Omega)$ is a local minimiser of \eqref{eqn distance functional} then it is also a local minimiser of \eqref{eqn distance functional relaxed scalar}, and so we may assume that $K$ is convex to prove Theorem \ref{Theorem scalar case regularity non convex distance functional}. We also mention that if $n=2$, Theorem \ref{Theorem scalar case regularity non convex distance functional} can be proven by duality methods and a result of Fonseca--Fusco, see \cite{Fonseca1997, Fonseca2002}.

To obtain a result in the general scalar case, we propose a Moser iteration scheme, combined with approximation methods based on methods of Fonseca--Fusco. To begin with, we state a general estimate which is a variant of the standard one, see \cite[Chapter 8]{Giusti}.
\begin{prop}\label{theorem fonseca fusco caccioppoli weighted}
Let $1<p<\infty$, and assume that $F,G \in C^2(\R^n)$ satisfy for some $\nu,L,\lambda,C_0>0$,
\begin{enumerate}[(i)]
    \item $\nu (1+|z|^2)^\frac{p-2}{2}I \leq F''(z) \leq \nu^{-1}(1+|z|^2)^{\frac{p-2}{2}}I$
    \item $0 \leq G(z) \leq L(1+|z|^2)$
    \item $\lambda I \leq G''(z) \leq LI$
    \item $G'(z) \otimes G'(z) \leq C_0G(z)G''(z).$
\end{enumerate}
Suppose that $u \in \WW^{1,p}_\text{loc}(\Omega)$ is a local minimiser of $\int_\Omega F(Dv) \dd x$. Then for every locally Lipschitz non-negative increasing function $h: [0,\infty) \rightarrow [0,\infty)$, and for every $\eta \in C_c^\infty(\Omega)$, $u$ satisfies the weighted Caccioppoli inequality
$$\lambda \int \eta^2 h(G(Du)) \sum_iF''(Du)[\p_iDu,\p_iDu] \dd x \lesssim C_0 \int h(G(Du)) G(Du)F''(Du)[D\eta, D\eta] \dd x$$
where the implied constant is independent of $\nu,L,\lambda$.
\end{prop}
\begin{proof}
From $(i)$, we know by a standard difference quotient argument \cite[Section 8]{Giusti} that $u \in \WW^{2, \min(2,p)}_\text{loc}(\Omega)$ with $(1+|Du|^2)^{\frac{p-2}{2}}|D^2u|^2 \in \LL^1_\text{loc}(\Omega).$ In particular this allows us to differentiate the Euler--Lagrange equation, so that for every $s \in \left\{1, \dots, n \right\}$ we have
$$\int_\Omega F''(Du) \p_s Du\cdot D\varphi \dd x =0$$
for any compactly supported $\varphi$ for which the integral makes sense. We define the truncation $h_T(t)=\min(h(t),T)$, and for each $s$ we use the test function $$\varphi_s=\eta^2 h_T(G(Du)) (\p_sG)(Du).$$ 
One can check that this is a valid test function using $(i),(ii),(iii)$ combined with the integrability of $D^2u$ mentioned above. To simplify our notation, we write $G=G(Du)$ and $G_s=(\p_sG)(Du)$, and similarly for $G',G''$. We compute
\begin{align*}
F''(Du)[\p_s Du, D\varphi_s] & = 2\eta h_T(G) F''(Du)[G_s\p_s Du, D\eta]\\
& +\eta^2 h_T'(G)F''(Du)[G_s\p_sDu, D^2uG']\\
& + \eta^2 h_T(G)F''(Du)[\p_s Du, D^2u G_s'].
\end{align*}
Noting that $\sum_s G_s\p_sDu = D^2 uG'$, and that $h_T'(G)F''(Du)[D^2uG',D^2uG'] \geq 0$ by convexity of $F$, after substituting this into the differentiated Euler--Lagrange equation and summing we have
$$\int \eta^2 h_T(G) \sum_s F''(Du)[\p_s Du, D^2u G_s'] \dd x \leq 2\int \eta h_T(G) \left|F''(Du)[D^2uG',D\eta]\right| \dd x.$$
Note that we can equivalently write $$ \sum_s F''(Du)[\p_s Du, D^2u G_s'] =\mathrm{Tr}(D^2uF''(Du)D^2uG'').$$
Once again using the convexity of $F$, we may apply the Cauchy--Schwarz inequality to the right-hand side to estimate
$$\left|F''(Du)[D^2uG',D\eta] \right| \leq \left(F''(Du)[D^2uG',D^2uG'] \right)^\frac{1}{2}\left(F''(Du)[D\eta,D\eta]\right)^\frac{1}{2}.$$
Using the fact that trace preserves inequalities between positive semidefinite matrices, our assumption $(iv)$ lets us estimate
\begin{align*}
F''(Du)[D^2uG',D^2uG']&=\mathrm{Tr}(D^2uF''(Du)D^2u (G'\otimes G'))\\
& \leq C_0G \mathrm{Tr}(D^2uF''(Du)D^2u G'').
\end{align*}
Since we now have $\mathrm{Tr}(D^2u F''(Du)D^2uG'')$ on both the left and right hand side of our inequality, we can apply Young's inequality and subtract in the usual fashion to get 
$$\int \eta^2 h_T(G) \mathrm{Tr}(D^2uF''(Du)D^2uG'') \dd x \lesssim4 C_0\int h_T(G)GF''(Du)[D\eta, D\eta] \dd x.$$
Finally we use $G'' \geq \lambda I$ to deduce 
$$\int \lambda \eta^2 h_T(G) \sum_i F''(Du)[\p_i Du, \p_i Du] \dd x \leq 4C_0 \int h_T(G)GF''(Du)[D\eta, D\eta] \dd x.$$
Now taking $T \rightarrow \infty$ and using the monotone convergence theorem, we obtain the claimed inequality.
\end{proof}
We wish to apply this to the distance functional. Recall that at points of differentiability of $\pi$ outside $K$, we have
 $$(p-1)|q|^{p-4}q \otimes q \leq \left(\frac{1}{p}F_K^p\right)''(z)\leq \max(p-1,1)|q|^{p-2}I.$$
 Considering this estimate, in order to carry out the Moser iteration starting from Proposition \ref{theorem fonseca fusco caccioppoli weighted}, it would be natural to take $F(z)=\frac{1}{p}F_K^p(z)$ and $G(z)=|q(z)|^2$. The issue is that neither of these are $C^2$ or satisfy the uniform convexity conditions. Therefore, we instead consider their regularisations. 
We define
the regularised distance 
$$d_\ep(z)=(\ep^2+(\rho_\ep *F_K^2)(z))^{1/2}$$
and then the functions
$$G_{\ep,\lambda}(z)=d_\ep^2(z)+\lambda|z|^2$$
$$F_\ep(z)=\frac{1}{p}d_\ep(z)^p+\frac{\ep}{p}(1+|z|^2)^{p/2}$$
where $\lambda, \ep \in (0,1)$. Clearly $d_\ep, G_{\ep,\lambda}, F_\ep$ are all smooth. The reason for mollifying $F_K^2$ instead of $F_K$ is not fundamental, but appears to make the algebra easier. The reason for adding $\ep^2$ in the distance $d_\ep$ is both to ensure it is smooth, and also for technical reasons that it is convenient to have $G(z)>0$ in the iteration. 
 \begin{lemma}\label{lemma technical estimates on regularisations}
Define $d_\ep$, $G_{\ep\,\lambda}, F_\ep$ as above for $0 < \ep,\lambda \leq 1$.
\begin{enumerate}[(i)]
\item $d_\ep \in C^\infty(\R^n)$ is convex, $d_\ep \geq \ep$,  $|d_\ep'(z)| \leq 1,$ and $0 \leq d_\ep(z)d_\ep''(z) \leq I$
    \item $\ep^2 \leq G_{\ep,\lambda}(z) \leq L(1+|z|^2)$
    \item $2\lambda I \leq G_{\ep,\lambda}''(z) \leq 2(1+\lambda)I$
    \item $G_{\ep,\lambda}'(z) \otimes G_{\ep,\lambda}'(z) \leq 2 G_{\ep,\lambda}(z)G''_{\ep, \lambda}(z)$
    \item $L^{-1}\ep (1+|z|^2)^{\frac{p-2}{2}}I +(p-1)d_\ep^{p-2}(z) d_\ep'(z) \otimes d_\ep'(z)\leq F_\ep''(z) \leq L \left[ \ep(1+|z|^2)^{\frac{p-2}{2}}I+d_\ep^{p-2}(z)\right]I$.
\end{enumerate}
where $L>0$ is independent of $\ep$ and $\lambda$.
\end{lemma}
\begin{proof}
The properties are all easy to verify but we give the details. Firstly note the convexity of $d_\ep$ follows from convexity of $F_K$. Recall that $0 \leq (F_K^2)^{''}(z) \leq 2I$, so by properties of mollification we also have $ 0 \leq (d_\ep^2)'' \leq 2I$. However we also have $(\frac{1}{2}d_\ep^2)''=d_\ep d_\ep''+d_\ep' \otimes d_\ep'$. Since both terms are positive semidefinite the estimates $(i)$ readily follow. The growth estimate $(ii)$ is immediate. Since $G_{\ep,\lambda}''(z)=(d_\ep^2)''+2\lambda I$, the estimates $(iii)$ also readily follow. Finally $(iv)$ follows from a similar Cauchy--Schwarz argument as in Lemma \ref{lemma inequality of q, pi and convexity of distance function}. To see $(v)$, we compute
$$F''_\ep(z)=d_\ep^{p-2}(z)((p-1)d_\ep'(z)\otimes d_\ep'(z) +d_\ep(z)d_\ep''(z))+\frac{\ep}{p}\left((1+|z|^2)^{p/2} \right)''.$$
The desired estimate follows by using $(i)$. 
\end{proof}
In particular, both $F_\ep$ and $G_{\ep,\lambda}$ satisfy the conditions of Proposition \ref{theorem fonseca fusco caccioppoli weighted}.

\begin{prop}\label{prop Lipschitz estimates for regularised problen scalar}
Let $1<p<\infty$, and assume that $v \in \WW^{1,p}_\text{loc}(\Omega)$ is a local minimiser of the regularised integral $\int_\Omega F_{\ep}(Dw) \dd x.$ Then $v \in \WW^{1,\infty}_\text{loc}(\Omega)$ and for all $B_R \Subset \Omega$ we have
\begin{equation}\label{eqn lipschitz estimate moser iteration}
\operatorname*{ess\,sup}_{B_{R/2}}|Dv|^p \leq C\left(1+\strokedint_{B_R}|Dv|^p \dd x \right)\end{equation}
where $C$ is independent of $\ep\in(0,\frac{1}{2})$.
\end{prop}
\begin{proof}
To simplify notation we write $d_\ep=d_\ep(Dv)$ and $G_{\ep,\lambda}=G_{\ep,\lambda}(Dv)$, so that $DG_{\ep,\lambda}=D(G_{\ep,\lambda}(Dv))$ for instance. For any locally Lipschitz non-negative increasing function $h:[0,\infty) \rightarrow [0,\infty)$, and every $\eta \in C_c^\infty(\Omega)$, we have by Proposition \ref{theorem fonseca fusco caccioppoli weighted} the estimate $$\lambda \int \eta^2 h(G_{\ep,\lambda}) \sum_iF_\ep''(Dv)[\p_iDv,\p_iDv] \dd x \lesssim C_0 \int h(G_{\ep,\lambda}) G_{\ep,\lambda}F_\ep''(Dv)[D\eta, D\eta] \dd x.$$ Next we apply $(v)$ in Lemma \ref{lemma technical estimates on regularisations}, using $Dd_\ep=d_\ep'(Dv)D^2v$, to get
$$\lambda \int \eta^2 h(G_{\ep,\lambda})d_\ep^{p-2}|Dd_\ep|^2 \dd x \lesssim \int h(G_{\ep,\lambda})G_{\ep,\lambda} \left[\ep(1+|Dv|^2)^{\frac{p-2}{2}}+d_\ep^{p-2} \right]|D\eta|^2 \dd x.$$
We wish to simplify this estimate by writing it all in terms of $d_\ep$. To start with, note that there exists $M>0$ depending only on $K$ such that $d_\ep(z) \approx |z| \approx (1+|z|^2)^{\frac{1}{2}}$ for all $|z| \geq M$, where the implied constants are all universal. We therefore are motivated to choose $h$ in such a way that $|Dv| \geq M$ on the support of $h(G_{\ep,\lambda}).$ To carry out the iteration, we also want $h(t) \approx t^m$, at least when $t$ is sufficiently large. We therefore define $h$ in the following way. Firstly, there exists $A>0$ depending only on $K$ such that $G_{\ep,\lambda}(z) \geq A$ implies $|z| \geq M$. We take an increasing locally Lipschitz  function $\chi \in C^\infty(0,\infty)$ so that $\chi=0$ on $[0,A]$ and $\chi=1$ on $[2A,\infty)$. For $m \geq 0$ we set $h(t)=\chi(t)t^m$. By construction, $|Dv| \geq M$ on the support of $h(G_{\ep,\lambda})$, so that we get
$$\lambda \int \eta^2 \chi(G_{\ep,\lambda})G_{\ep,\lambda}^m d_\ep^{p-2}|Dd_\ep|^2 \dd x \lesssim \int \chi(G_{\ep,\lambda})G_{\ep,\lambda}^{m+1}d_\ep^{p-2}|D\eta|^2 \dd x.$$
Next, we want to choose $\lambda$ so that $G_{\ep,\lambda}^m$ is comparable to $d_\ep^{2m}$ whenever $|Dv| \geq M$. Recalling that $G_{\ep,\lambda}(z)=d_\ep^2(z)+\lambda|z|^2$ and $d_\ep \approx |Dv|$ if $|Dv| \geq M$, we clearly have
$$d_\ep^2 \leq G_{\ep,\lambda} \leq (1+C\lambda)d_\ep^2. $$
Taking for instance $\lambda_m=\frac{1}{m+1}$, it follows that $G_{\ep,\lambda_m}^m \approx d_{\ep}^{2m}$ on the support of $\chi(G_{\ep,\lambda}).$ It follows that
$$\int \eta^2 \chi(G_{\ep,\lambda_m})d_\ep^{2m+p-2}|Dd_\ep|^2 \dd x \lesssim (m+1)\int d_\ep^{2m+p}|D\eta|^2 \dd x.$$
This is very nearly the desired inequality to carry out Moser iteration, except for the presence of the cut-off $\chi(G_{\ep,\lambda_m})$. However, this can be easily fixed by defining a truncation, for instance $$T_\ep=\max(d_\ep, \sqrt{2A}).$$ If $d_\ep \leq \sqrt{2A}$ then $DT=0$, while if $d_\ep \geq \sqrt{2A}$ then $G_{\ep,\lambda_m} \geq 2A$ so that $\chi(G_{\ep,\lambda_m})=1$. Consequently, 
$$\int \eta^2 T_\ep^{2m+p-2}|DT_\ep|^2 \lesssim (m+1)\int T_\ep^{2m+p}|D\eta|^2 \dd x.$$
From here the argument is standard, see \cite[Proposition 8.1]{Giusti}. Indeed, we may rewrite the above inequality as 
$$\int \eta^2 |D(T_\ep^k)|^2 \dd x \lesssim k^3 \int T_\ep^{2k}|D\eta|^2 \dd x$$
where $k=m+\frac{p}{2}$.
One combines this with the Sobolev--Poincaré inequality applied to $\eta T^k$, and after iterating this estimate on concentric balls with corresponding choices for $\eta$, one obtains
$$\sup_{B_{R/2}} T_\ep \leq C\left( \strokedint_{B_R}T_\ep^p \dd x\right)^\frac{1}{p}$$
where $C$ is a universal constant. 
Unpacking the definitions we see that $T_\ep \approx 1+|Dv|$, so this estimate implies \eqref{eqn lipschitz estimate moser iteration}.
\end{proof}
The last step is to pass the regularity from the regularised problem to the actual problem. The standard method one employs is to consider minimisers $u_\ep \in u+\WW_0^{1,p}$ of the regularised problem. From the previous proposition we know $u_\ep$ are Lipschitz with an $\ep$-independent estimate. Since $F_\ep \rightarrow \frac{1}{p}F_K^p$ locally uniformly, this is enough to carry the estimate over to $u$. An interesting subtlety here is that, due to lack of strict convexity, we can't necessarily conclude that $Du_\ep$ converges to $Du$.

\begin{proof}[Proof of Theorem \ref{Theorem scalar case regularity non convex distance functional}]
Let $u_\ep \in u+W_0^{1,p}(B_R)$ be a minimiser of $\int_{B_R}F_\ep(Dw) \dd x$ over this Dirichlet class. Let us record that $$\frac{1}{p}F_K^p(z) \leq F_\ep(z)$$ using Jensen's inequality, and clearly $$|z|^p \lesssim 1+F_K^p(z) \lesssim 1+F_\ep(z) \lesssim 1+|z|^p$$ Using this and  minimality of $u_\ep$ we easily get
$$\strokedint_{B_R}|Du_\ep|^p \dd x \lesssim \strokedint_{B_R}1+|Du|^p \dd x$$
and hence Proposition \ref{prop Lipschitz estimates for regularised problen scalar} also gives
$$\sup_{B_{R/2}}|Du_\ep|^p \lesssim \strokedint_{B_R}1+|Du|^p \dd x.$$
Using weak compactness of $\WW^{1,p}$ and weak$^*$ compactness of $\WW^{1,\infty}$, we know there exists $v \in u+W_0^{1,p}(B_R)$ such that $u_\ep \rightarrow v$ weakly in $\WW^{1,p}(B_R)$ and where $$\sup_{B_{R/2}}|Dv|^p \lesssim \strokedint_{B_R}1+|Du|^p \dd x.$$ Next, we see that $v$ is also a minimiser of $\int_{B_R}F_K^p(Dw) \dd x$ over its Dirichlet class; indeed using weak lower semicontinuity of the integrand, local uniform convergence of $F_\ep \rightarrow \frac{1}{p}F_K^p$, minimality of $u_\ep$, and $\frac{1}{p}F_K^p \leq F_\ep$, we have
\begin{align*}
\frac{1}{p}\int_{B_R}F_K^p(Dv) \dd x &\leq \liminf_{\ep \searrow 0} \frac{1}{p}\int_{B_R}F_K^p(Du_\ep) \dd x \\
& \leq \liminf_{\ep \searrow 0}\int_{B_R}F_\ep(Du_\ep) \dd x\\
& \leq \liminf_{\ep \searrow 0}\int_{B_R}F_\ep(Du) \dd x\\
&=\frac{1}{p}\int_{B_R}F_K^p(Du) \dd x.
\end{align*}
Consequently, both $u$ and $v$ are minimisers of $\int_{B_R}F_K^p(Dw)\dd x$ over the same Dirichlet class. This means $q(Du)=q(Dv)$ almost everywhere. Indeed, 
$$\int_{B_R}F_K^p\left(\frac{Du+Dv}{2}\right) \dd x \leq \int_{B_R}\left|\frac{q(Du)+q(Dv)}{2} \right|^p \leq \int_{B_R}\frac{F_K^p(Du)+F_K^p(Dv)}{2} \dd x$$
and minimality forces this to be an equality, In particular $\left|\frac{q(Du)+q(Dv)}{2} \right|^p=\frac{|q(Du)|^p+|q(Dv)|^p}{2}$ almost everywhere which implies $q(Du)=q(Dv)$. Since $q(z)=z-\pi_K(z)$ and $\pi_K$ is bounded, clearly this implies the boundedness of $Du$ together with the required estimate.
\end{proof}
Theorem \ref{Theorem scalar case regularity non convex distance  functional} is sharp. Recall that if $K = \left\{a\right\}$ is a singleton, then minimisers are, up to affine perturbations, just $p$-harmonic maps and therefore locally $C^{1,\alpha}$. From Proposition \ref{prop rank one connections mean discontinuous derivatives}, it follows in the scalar case that whenever $K \subset \R^n$ contains two distinct points, we do not expect $C^1$ regularity of minimisers. To conclude this section, we remark the following situation where the vectorial problem decouples into scalar problems.
\begin{prop}\label{prop vectorial regularity decouples}
Let $K = K_1 \times \dots \times K_N \subset \R^{N \times n}$ where each $K_\alpha \subset \R^n$ is compact and convex. Suppose that $u \in \WW^{1,2}_\text{loc}(\Omega;\R^N)$ is a local minimiser of $\int_{\Omega}F_K^2(Dv) \dd x.$ Then $u \in \WW_\text{loc}^{1,\infty}(\Omega;\R^N)$ together with the estimate
$$\sup_{B_{R/2}}|Du|^2 \ \leq C \left(1+\strokedint_{B_R}|Du|^2 \dd x \right)$$
for all $B_R \Subset \Omega.$
\end{prop}
\begin{proof}
Simply observe that $F_K^2(Du) = \sum_{\alpha=1}^N\dist^2(Du^\alpha,K_\alpha)$, so that each $u^\alpha$ minimises $\int_\Omega \dist^2(Dv^\alpha,K_\alpha) \dd x.$ By Theorem \ref{Theorem scalar case regularity non convex distance functional} it follows $u^\alpha$ is locally Lipschitz for all $\alpha$.
\end{proof}

\section{$\BMO$ regularity in the non-convex vectorial case}
In the rest of this paper we consider the vectorial case, namely where $K \subset \R^{N \times n}$ and $N>1$. The regularity of minimisers of $\int_\Omega F_K^p(Du) \dd x $ is more subtle, and as we will see the geometry of $K$ plays an important role. The main result in the general case of compact sets $K \subset \R^{N \times n}$ is Theorem \ref{Theorem BMO regularity for distance functional}, which we recall says that minimisers of the relaxation $\int_\Omega (F_K^p)^{qc}(Dv) \dd x $ have locally $\BMO$ gradients. The proof of this follows immediately from the estimate
$$|\nab (F_K^p)^{qc}(z)-pA_p(z)| \leq L(1+|z|)^{p-2}$$
and the following result which is proved in \cite{Strong2026}.
\begin{prop}
    Let $1<p<\infty$, and suppose $u \in \WW^{1,p}(\Omega;\R^N)$ is a local minimiser of $\int_\Omega F(Dv) \dd x$ where $F \in C^1(\R^{N \times n})$ is a continuous integrand satisfying $$\limsup_{z \rightarrow \infty}\frac{|F'(z)-pA_p(z)|}{|z|^{p-2}}<\infty.$$ Then $$Du \in \BMO_\text{loc}(\Omega;\R^{N \times n}).$$
\end{prop}
The sharpness of Theorem \ref{Theorem BMO regularity for distance functional} is already given in \cite{DolzmannKristensenZhang2013}, where non--Lipschitz minimisers for several functionals are considered. For example, non--Lipschitz minimisers of the quadratic distance functional are constructed in the case $K=SO(2) \subset \R^{2 \times 2}$.

\section{Sobolev regularity in the convex vectorial case}
For the rest of this paper we take $K \subset \R^{N \times n}$ to be a compact convex set, where $N>1$. To begin with, we demonstrate via a difference quotient argument that the stress tensor has meaningful derivatives. Recall that for $1<p<\infty$, if $u \in \WW^{1,p}_\text{loc}(\Omega;\R^N)$ is a local minimiser of $\int_\Omega F_K^p(Dv) \dd x$ then $u$ satisfies the Euler--Lagrange equation 
\begin{equation}\label{eqn euler lagrange equation section 10}
-\Div \sigma=0
\end{equation}
where $$\sigma=\sigma(Du)=A_p(q(Du))$$ denotes the stress. By convexity of $F_K$, any solution to the Euler--Lagrange equation is automatically a minimiser in its Dirichlet class. To derive difference quotient estimates, we will use the monotonicity estimate
$$\Delta_j^h\sigma \cdot \Delta_j^h \pi(Du)\geq 0$$
where $\Delta_j^hf(x)=f(x+he_j)-f(x)$. Note that this estimate follows immediately from Lemma \ref{lemma monotonicity estimates for sigma}. We recall the Niksolkii--Besov space $N^{s,p}(\Omega)=B^{s}_{p, \infty}(\Omega)$ is defined as the set of functions $u \in \LL^p(\Omega)$ such that
$$[u]_{N^{s,p}(\Omega)}^p
:=\sum_{j=1}^n\sup_{|h|<1}
\int_{\Omega}
\frac{|\Delta_j^h u|^p}{|h|^{ps}} \dd x
<\infty
$$
where we make the abuse of notation that $\Delta_j^h u(x)$ equal zero if $x+he_j \notin \Omega$.

\begin{theorem}\label{theorem difference quotient superquadratic distance function}
Let $1<p<\infty$ and $K \subset \R^{N \times n}$ compact and convex. Suppose that $u \in \WW^{1,p}_\text{loc}(\Omega;\R^N)$ is a solution of the Euler--Lagrange equation \eqref{eqn euler lagrange equation section 10}.
\begin{enumerate}[(i)]
    \item If $p \geq 2$ then $V_{p'}(\sigma) \in \WW^{1,2}_\text{loc}(\Omega;\R^{N \times n})$ together with the estimate 
\begin{equation}\label{eqn difference quotients p>=2}
\int_{B_{R/2}}|DV_{p'}(\sigma)|^2 \dd x \leq \frac{C}{R^2}\int_{B_R}|q(Du)|^{p-2}|Du-z_0|^2\dd x
\end{equation}
for all $B_R \Subset \Omega$ and $z_0 \in \R^{N \times n}$.
\item If $1<p<2$ then $V_{p'}(\sigma) \in \NN_\text{loc}^{\frac{p}{2},2}(\Omega;\R^{N \times n} )$ together with the estimate
\begin{equation}\label{eqn difference quotients 1<p<2}
[V_{p'}(\sigma)]_{\NN^{\frac{p}{2},2}(B_{R/2})}^2 \leq \frac{C}{R^p}\int_{B_R}|Du-z_0|^p \dd x\end{equation}
for all $B_R \Subset \Omega$ and $z_0 \in \R^{N \times n}$.
\end{enumerate}

\end{theorem}
\begin{proof}
We fix $\eta \in C_c^\infty(\Omega;[0,1])$ and take $\Delta_j^{-h}(\eta^2\Delta_j^h(u-l))$ as a test function in \eqref{eqn euler lagrange equation section 10}, where $l$ is an affine map with $Dl=z_0$. Expanding this out, this gives
$$\int \eta^2 \Delta_j^h \sigma \cdot \Delta_j^h Du \dd x = -2\int \eta \Delta_j^h \sigma \cdot  \Delta_j^h(u-l) \otimes D\eta  \dd x.$$
Since $Du=A_{p'}(\sigma)+\pi(Du)$, the monotonicity estimate $\Delta_j^h\sigma \cdot \Delta_j^h\pi_K(Du) \geq 0$ gives
$$\int \eta^2 \Delta_j^h \sigma \cdot \Delta_j^hA_{p'}(\sigma) \dd x \lesssim \int_{\Omega } \eta |\Delta_j^h \sigma||D\eta ||\Delta_j^h(u-l)| \dd x.$$
By Lemma \ref{lemma equivalences of auxiliary V functinos} the equivalence $\Delta_j^h \sigma \cdot \Delta_j^h A_{p'}(\sigma) \approx |\Delta_j^hV_{p'}(\sigma)|^2$ holds and therefore $$\int \eta^2 |\Delta_j^hV_{p'}(\sigma)|^2 \dd x \lesssim \int_{\Omega } \eta |\Delta_j^h \sigma||D\eta ||\Delta_j^h(u-l)| \dd x.$$
 From here we split between the superquadratic and subquadratic cases. \\
 \underline{Case 1: $p \geq 2$}.
 By Lemma \ref{lemma equivalences of auxiliary V functinos} we note the equivalence $$|A_p(z)-A_p(w)| \approx (|z|+|w|)^{\frac{p-2}{2}}|V_p(z)-V_p(w)|.$$ Unpacking the definition $\sigma =A_p(q(Du))$ and using $V_p(q(Du))=V_{p'}(\sigma)$, this yields $$|\Delta_j^h \sigma | \approx (|q(Du(x+he_j))|+|q(Du(x)|)^{\frac{p-2}{2}} |\Delta_j^h V_{p'}(\sigma)|.$$
Applying Young's inequality and subtracting from the right-hand side, we arrive at 
$$\int \eta^2 |\Delta_j^h V_{p'}(\sigma)|^2 \dd x \lesssim \int |D\eta|^2 (|q(Du(x+he_j))|+|q(Du(x))|)^{p-2}|\Delta_j^h(u-l)|^2 \dd x.$$
At this point we simply choose $\eta = 1$ on $B_{R/2}$ with $\eta=0$ outside $B_R$ and $|D\eta| \lesssim \frac{1}{R}$. It remains to justify taking the limit $h \searrow 0$, which is justified here because $p \geq 2$. Indeed, by the difference quotient lemma it follows $\frac{|\Delta_j^h(u-l)|^2}{h^2} \rightarrow |\p_j(u-l)|^2$ strongly in $\LL^{\frac{p}{2}}$, and similarly because $q$ is Lipschitz and translations converge in $\LL^p$, we have $(|q(Du(x+he_j))|+|q(Du(x))|)^{p-2}\rightarrow(2q(Du))^{p-2}$ in $\LL^{\frac{p}{p-2}}$. After dividing by $h^2$, it follows by Hölder's inequality the right-hand side converges, and therefore the left-hand side is uniformly bounded. This entails weak convergence of $\frac{\Delta_j^hV_{p'}(\sigma)}{h}$ in $\LL^2$, and therefore by weak lower semicontinuity we obtain \eqref{eqn difference quotients p>=2}.\\
\underline{Case 2: $1<p<2.$}
This time we write $\eta = \eta^{\frac{2}{p'}}\cdot \eta^{\frac{2-p}{p}}$ and apply Young's inequality with exponent $p'$ to get 
$$\int \eta |\Delta_j^h \sigma||D\eta||\Delta_j^h(u-l)| \dd x \lesssim \ep \int \eta^2 |\Delta_j^h\sigma|^{p'}\dd x+\int \eta^{2-p}|D\eta|^p |\Delta_j^h(u-l)|^p \dd x.$$
Since $p<2$ we also have
$$|\Delta_j^h V_{p'}(\sigma)|^2 \gtrsim |\Delta_j^h \sigma|^{p'}.$$ Using this, we may subtract from both sides as usual to get
$$\int \eta^2 |\Delta_j^h V_{p'}(\sigma)|^2 \dd x \lesssim \int \eta^{2-p}|D\eta|^p|\Delta_j^h(u-l)|^p \dd x.$$
Choosing $\eta$ as in the previous proof and applying the difference quotient lemma, this gives
$$\int_{B_{R/2}}|\Delta_j^hV_{p'}(\sigma)|^2 \dd x \leq C\frac{|h|^p}{R^p}\int_{B_R}|Du-z_0|^p dd x$$
from which the claim follows.
\end{proof}
We remark that \cite{DolzmannKristensenZhang2013} carried out the above difference quotient argument in the case $p=2$. Note that we can equivalently write \eqref{eqn difference quotients p>=2} as $$\int_{B_{R/2}}|DV_{p'}(\sigma)|^2 \dd x \leq \frac{C}{R^2}\int_{B_R}|q(z_0)|^{p-2}|Du-z_0|^2+|Du-z_0|^p\dd x$$
since $q$ is Lipschitz. Because $Du \in \BMO_\text{loc}$, this estimate entails a Morrey type estimate for $DV_{p'}(\sigma)$. Indeed, fixing some ball $B_{R_0}$, it follows
$$\int_{B_{R/2}}|DV_{p'}(\sigma)|^2 \dd x \leq CR^{n-2}(|q((Du)_R)|^{p-2}+1))$$
where $C$ depends on $[Du]_{\BMO(B_{R_0})}$. Since $q$ is Lipschitz and $\BMO$ functions have comparable means on the balls $B_R$ and $B_{2R}$, it follows that $$|q((Du)_R)| \lesssim 1+\log \frac{1}{R}$$ which therefore yields a logarithmic Morrey estimate
$$\int_{B_{R/2}}|DV_{p'}(\sigma)|^2 \dd x \lesssim C R^{n-2} \left(1+\log \frac{1}{R}\right)^{p-2}$$
which holds for $p \geq 2$.

\section{Counterexample to Lipschitz regularity in the vectorial convex case}
In this section we prove Theorem \ref{theorem counterexample lipschitz vectorial first one}, namely we show that Theorem \ref{Theorem scalar case regularity non convex distance functional} does not generalise to the vectorial case $N >1$. Our construction relies on constructing a compact set $K \subset \R^{N \times n}$ and a non-Lipschitz map $u \in \WW^{1,2}$ such that:
\begin{enumerate}[(i)]
    \item $Du=P+\nu$ where $\nu$ is divergence free and $P \in K$.
    \item $\nu(x) \cdot (k-P(x)) \leq 0$ for all $k \in K$ and $x \in \R^n$.
\end{enumerate}
If these two conditions hold then $u$ is a local minimiser of the square distance functional. Indeed, condition $(ii)$ says that $\nu \in N_K(P)$. Because $N_K(P)=\pi^{-1}(P)-P$, $(ii)$ implies $P=\pi(P+\nu)$. Then $(i)$ gives $\pi(Du)=P$ and $\nu=\sigma$ where $\sigma=Du-\pi(Du)$. Since $\nu$ is divergence free, the Euler--Lagrange equation $-\Div \sigma=0$ holds and therefore $u$ is a local minimiser of $\int_\Omega F_K^2(Dv) \dd x$. 
A natural first idea would be to consider a map such as $$v(x)=(n-1) x \log \frac{1}{|x|}$$
on $B_1(0) \subset \R^n$.
Writing $r=|x|$ and $\theta=\frac{x}{r}$, one computes
$$Dv=-(n-1)\theta \otimes \theta + (n-1)\log \frac{1}{r}I.$$
We further compute $$-\Div \left(\log \frac{1}{r}I\right)=D \left(\log \frac{1}{r} \right) = -\frac{\theta}{r}$$
$$\Div(\theta \otimes \theta)=(D\theta)\theta+\theta \Div \theta=\frac{n-1}{r}\theta.$$
Inspecting this, a natural decomposition is therefore
$$Dv=P_v+\nu_v$$
where
$$P_v(x)=-n \theta \otimes \theta$$
$$\nu_v(x)=(n-1)\log \frac{1}{r}I+\theta \otimes \theta.$$
We see that $\nu_v$ is divergence free. Note that $P_v$ is $0$-homogeneous so we may think of it as a function of $\theta$. Since we want $P_v \in K$, we could consider taking $$K=\mathrm{conv.hull}(P_v(\theta):\theta \in S^{n-1})$$ so that condition $(i)$ is automatic. Unfortunately, condition $(ii)$ then fails; indeed one computes
$$\nu_v(x)\cdot (P_v(\eta)-P_v(\theta))=n(1-(\theta \cdot \eta)^2) \geq 0$$ for all 
$\theta, \eta \in S^{n-1}$, where $x=\theta r$. 
It turns out that one can fix this by considering $$u(x)=(v(x),w(x))$$ where $w$ is a suitably chosen map so that the supporting inequality is flipped.
\begin{proof}[Proof of Theorem \ref{theorem counterexample lipschitz vectorial first one}]
We take $v(x)=(n-1)x \log \frac{1}{r}$ and $P_v, \nu_v$ as above. We consider the $1$-homogeneous map
$$w(x)=\lambda \frac{x \otimes x}{|x|}=\lambda r \theta \otimes \theta$$
for some $\lambda \in \R$ chosen later. We compute the derivative $$Dw =\lambda \left[\frac{1}{r}D(x \otimes x) -\frac{(x \otimes x) \otimes x}{|x|^3} \right] = \lambda [D(x \otimes x)(\theta)-(\theta \otimes \theta) \otimes \theta]$$
where the notation $D(x \otimes x)(\theta)$ means $f'(\theta)$ with $f(x)=x \otimes x$. We need to find a convenient splitting of $Dw$ into a bounded and divergence free part. Of course, $Dw$ itself is bounded so there is a trivial splitting of $P_w=Dw$ and $\nu_w=0$. However this choice does not change the supporting inequality $(ii)$. Instead, we observe the algebraic identities 
$$\Div[D(x \otimes x)(\theta)]=\frac{2}{r}(I-\theta \otimes \theta)$$
$$\Div ((\theta \otimes \theta)\otimes \theta)=\frac{n-1}{r}\theta \otimes \theta$$
$$\Div(I \otimes \theta)=\frac{n-1}{r}I.$$
Therefore a natural divergence free candidate is 
$$\nu_w := D(x \otimes x)(\theta)-\frac{2}{n-1}I\otimes \theta + \frac{2}{n-1}(\theta \otimes \theta)\otimes \theta.$$
With this choice, and defining $P_w = Dw-\nu_w$, we set $$u=(v,w), P=(P_v,P_w), \nu=(\nu_v,\nu_w)$$
and $$K=\mathrm{conv.hull}(P(\theta):\theta \in S^{n-1}).$$
We now check the supporting inequality $(ii)$. 
A tedious but simple computation yields
$$\nu(x) \cdot (P(\eta)-P(\theta))=(1-c)\left[n(1+c)-2(n+1)(\lambda-1)+\left(\lambda+\frac{2}{n-1}\right)\left(2+\frac{2}{n-1}c(c+1)\right) \right]$$
where $x=r\theta$ and $c=\eta \cdot \theta$. 
For any $n \geq 3$, one can choose $\lambda$ sufficiently large so that the right-hand side is $\leq 0$ for all $-1\leq c \leq 1$. For example, if $n=3$ and $\lambda=\frac{9}{2}$ then the right-hand side becomes $-\frac{1}{2}(1-c)^2(11c+28) \leq 0.$ By definition of $K$, this shows that $(ii)$ holds, and hence $u$ is a solution of the Euler--Lagrange equation. For $n=3$, we may view $u : B_1(0) \rightarrow \R^{3+6}$ since $w$ is symmetric.
\end{proof}
The above construction provides a counterexample for all $n \geq 3$, but not $n=2$. This is because the map $u(x)=x \log \frac{1}{|x|}$ does not have a $\WW^{1,2}$ stress, which is required by Theorem \ref{theorem difference quotient superquadratic distance function}. To construct a counterexample for $n=2$, one should instead consider a map where $u(x) \sim x \log \log \frac{1}{|x|}$ for small $|x|$. With the assistance of AI tools, the reader can be convinced that several planar counterexamples are possible. However, we have been unable to find a suitably simple construction, so we intend to report on this elsewhere.
\section{Lipschitz regularity assuming a supporting ball condition}
In this section we prove Theorem \ref{theorem lipschitz regularity for distance functional under supporting ball condition}. Our proof relies on the following general result from \cite{Strong2026}. For $F \in C^1(\R^{N \times n})$ we define
$$R_a(z)=\left[F(z)-F(a)-F'(a)\cdot (z-a)\right]-\left[|z|^p-|a|^p-pA_p(a)\cdot (z-a)\right].$$
\begin{prop}\label{theorem lipschitz regularity under approximate C^2 sense}
Suppose that $u \in \WW^{1,p}(\Omega;\R^N)$ is a local minimiser of $\int_\Omega F(Dv) \dd x$ where $F \in C^1(\R^{N \times n})$ satisfies \begin{equation}\label{eqn approximate dirichlet in C^2 sense with taylor}
\lim_{|a| \rightarrow \infty} \sup_{z \in \R^{N \times n}}\frac{|R_a(z)|}{|V_p(z)-V_p(a)|^2}=0.
\end{equation}
Then $u \in \WW^{1,\infty}_\text{loc}(\Omega;\R^N)$ and $$\sup_{B_{R/2}}|Du|^p \ \leq C \left(1+\strokedint_{B_R}|Du|^p \dd x \right)$$
for all $B_R \Subset \Omega.$
\end{prop}
Therefore, in order to show Lipschitz regularity of minimisers of the distance functional, we simply must verify that condition \eqref{eqn approximate dirichlet in C^2 sense with taylor} holds. When $K$ satisfies the uniform supporting ball condition, the estimates from Lemma \ref{lemma estimates SBC} will imply this condition.

\begin{proof}[Proof of Theorem \ref{theorem lipschitz regularity for distance functional under supporting ball condition}] 
We write $G(z)=F_K^p(z)-|z|^p$. By Lemma \ref{lem: Bregman divergence for V function}, the Theorem follows if we check that
$$|G(z)| \leq L(1+|z|^{p-1})$$
for all $z \in \R^{N \times n}$, and that
$$|G'(z)| \lesssim |z|^{p-2}, \qquad |G''(z)| \lesssim |z|^{p-3}$$
for $|z|$ sufficiently large. Indeed, if this holds then there exist $C,M>0$ such that $$\sup_{z \in \R^{N \times n}}\frac{|R_a(z)|}{|V_p(z)-V_p(a)|^2} \leq \frac{C}{|a|}$$ for all $|a| >M$, which obviously implies \eqref{eqn approximate dirichlet in C^2 sense with taylor}.  Recall that for $z$ outside of $K$ we know that $G''$ makes sense as a locally bounded measurable function. The first two estimates do not require uniform SBC, and just follow from Lemma \ref{Lemma growth conditions of distance functional}.
If $K$ is uniformly SBC with radius $R$, by Lemma \ref{lemma estimates SBC} we have almost everywhere the bound $$|\pi'(z)| \leq \frac{R}{R+|q(z)|}.$$ In particular, for large $|z|$ this implies $|\pi'(z)| \lesssim \frac{1}{|z|}$. Note that $$G''(z)=p\left[A_p'(q(z))-A_p'(z)-A_p'(q(z))\pi'(z)\right].$$ Using the fact that $|A_p'(z)| \lesssim |z|^{p-2}$ and $|A_p''(z)| \lesssim |z|^{p-3}$ for large $|z|$, as well as the boundedness of $\pi_K(z)$, the desired bound $$|G''(z)|\lesssim |z|^{p-3}$$
for large $|z|$ follows at once.
\end{proof}

\section{Further regularity results}
In this final section, we consider some geometric conditions on $K$ which allow us to deduce some further regularity for minimisers of the distance functional. For simplicity we mostly work in the quadratic case $p=2$, although it would be interesting to generalise the results to the full range $1<p<\infty$. 
\subsection{Rank-one connections}
The following basic observation tells us that the presence of rank-one connections in $K$ is an obstruction to regularity.
\begin{prop}\label{prop rank one connections mean discontinuous derivatives}
Let $1<p<\infty$ and suppose that $K \subset \R^{N \times n}$ has a rank-one connection. Then there exists a local minimiser $u \in \WW^{1,p}_\text{loc}(\Omega;\R^N)$ of $\int_\Omega F_K^p(Dv) \dd x$ such that $Du$ is not continuous.
\end{prop}
\begin{proof}
By translating we may assume $0 \in \Omega$. Let $A, B \in K$ be distinct matrices such that $B-A=a \otimes b.$ Constructing the usual laminate $u(x)=Ax+a(x\cdot b)_+$, we see that $Du=\begin{cases}A, & x\cdot b <0\\
B, & x\cdot b >0.
\end{cases}$
Then $Du$ is not continuous, but $F_K(Du)=0$ so $u$ is a local minimiser. 
\end{proof}
In the scalar case, lack of rank-one connections means that $K$ is a singleton, and therefore minimisers have $u \in C^{1,\alpha}_\mathrm{loc}$. In the theory of rigidity, it is well known  (see \cite{Rindler}) that if $L$ is a linear subspace containing no rank-one connections, then any exact solution of the differential inclusion $Du \in L$ is smooth. This is a consequence of the fact that the projection operator onto $L^\perp$ is suitably elliptic. When $K$ is convex, we can use this in our setting to carry out difference quotients and obtain Sobolev regularity on $u$. In the case $n=2$ this leads to continuity of the gradient. For simplicity we work with the case $p=2$. The following convex analysis fact is useful. Although it is probably well known to experts, we cannot find a precise reference so include a short proof here.

\begin{lemma}
Let $K \subset \R^{N \times n}$ be compact and convex containing $0$. Then the following are equivalent:
\begin{enumerate}[(i)]
    \item $K$ contains no rank-one connections
    \item There is a linear subspace $L \subset \R^{N \times n}$ such that $K \subset L$ and $L$ contains no rank-one matrices.
\end{enumerate}
\end{lemma}
\begin{proof}
Suppose that $(i)$ holds and let $L=\mathrm{span}  K$. Take a basis
$A_1,\ldots,A_d\in K$ of $L$. Since $0\in K$, it follows $\pm A_i \in K-K$ for all $i$. Moreover $K-K$ is convex, so that
$$\mathrm{conv}\{\pm A_1,\ldots,\pm A_d\}\subset K-K.$$
Thus $K-K$ contains a neighbourhood of $0$ relative
to $L$. If $L$ contained some rank-one matrix
$A$, then this means that $tA\in K-K$ for sufficiently small $t>0$. But this means $K$ has a rank-one connection, contradicting $(i)$. 
\end{proof}

\begin{prop}
Suppose that $K \subset \R^{N \times n}$ is compact and convex with no rank-one connections. Suppose that $u \in \WW^{1,2}_\text{loc}(\Omega;\R^N)$ is a local minimiser of $\int_\Omega F_K^2(Dv) \dd x.$ Then $u \in \WW^{2,2+\ep}_\text{loc}(\Omega;\R^N)$ for some $\ep>0$ depending on $K$ and universal constants. In particular, if $n=2$ then $u \in C^{1,\alpha}_\mathrm{loc}$ for some $\alpha>0$. 
\end{prop}

\begin{proof}
By translating we may assume $K \subset L$ where $L$ is a linear subspace with no rank-one matrices. Let $P$ be the orthogonal projection onto $L^\perp$. If $P(\lambda)=0$ then $\lambda \in L$ and therefore $\mathrm{rank }\lambda \neq 1$. Consequently $$c := \min_{|\lambda|=1, \rank \lambda=1}|P(\lambda)|^2>0$$
Equivalently, $|P(\lambda)|^2 \geq c|\lambda|^2$ for all $\rank \lambda=1$, which is precisely the Legendre-Hadamard condition associated to $P$ (noting $|P(z)|^2 = P(z) \cdot z$). Therefore, $P$ is quasielliptic in the sense that $$\int_{B_R} |P(D\varphi)|^2 \dd x \geq c \int_{B_R}|D\varphi|^2 \dd x $$ holds for all $\varphi \in W_0^{1,2}(B_R;\R^N).$ Following standard linear elliptic theory, to derive an estimate we take $\varphi = \eta \Delta_j^h(u-l)$ in the above inequality, where $\eta \in C_c^\infty(\Omega)$ and $l$ is an affine map. Using the linearity of $P$ and $|P(z)| \leq 
|z|$, this yields
$$\int \eta^2|\Delta_j^hDu|^2 \dd x \lesssim \int|D\eta|^2|\Delta_j^h(u-l)|^2+\eta^2|\Delta_j^hP(Du)|^2 \dd x.$$
Next, recall that from the proof of Theorem \ref{theorem difference quotient superquadratic distance function} we have the estimate
$$\int \eta^2 |\Delta_j^h q(Du)|^2 \dd x \lesssim \int |\Delta_j^h (u-l)|^2 |D \eta|^2 \dd x.$$
Since $P(Du)=P(q(Du))$, we have
$$|\Delta_j^h P(Du)|^2 = |P(\Delta_j^hq(Du))|^2 \leq |\Delta_j^h q(Du)|^2.$$
Putting this all together gives
$$\int \eta^2 |\Delta_j^h Du|^2 \dd x \lesssim \int |D\eta|^2 |\Delta_j^h(u-l)|^2 \dd x.$$
Suitably choosing $\eta$ and applying the difference quotient lemma gives
$$\int_{B_{R/2}}|D^2u|^2 \dd x \leq C\int_{B_R}\left|\frac{Du-z_0}{R} \right|^2\dd x$$
for all $z_0 \in \R^{N \times n}$, where $C$ is independent of $z_0$. By standard methods (see \cite{Giusti}), after selecting $z_0=(Du)_R$ and applying the Sobolev--Poincaré inequality and Gehring's lemma, we deduce that $u \in \WW^{2,2+\ep}_\mathrm{loc}(\Omega;\R^N)$ together with the estimates
$$\left( \strokedint_{B_{R/2}}|D^2u|^{2+\ep}\dd x\right)^{\frac{1}{2+\ep}} \lesssim \left(\strokedint_{B_{3R/4}}|D^2u|^2 \dd x\right)^{1/2} \lesssim \left(\strokedint_{B_R} \left|\frac{Du-(Du)_R}{R}\right|^2 \dd x\right)^{1/2}.$$
\end{proof}
If $K$ is not convex, it is possible for $K$ to have no rank-one connections and minimisers to be non-Lipschitz. Indeed, this happens if  $K=\left\{I,-I \right\}$ or $K=SO(2)$, as shown in  \cite{DolzmannKristensenZhang2013}. In the convex case with $n \geq 3$, $N \geq 2$, we believe it to be possible to construct non-Lipschitz minimisers for some convex $K$ with no rank-one connections, by a similar method to our proof of Theorem \ref{theorem counterexample lipschitz vectorial first one}. However, we intend to report on this elsewhere. 
\subsection{Line segments}
A line segment is a set of the form $$K=[A,B] = \left\{tA+(1-t)B: 0 \leq t \leq 1 \right\}.$$
Clearly we could also write a general line segment as $[B-A,B+A]$. Since $$\dist(Du, [B-A,B+A])=\dist(D(u-Bx),[-A,A]),$$ 
it suffices to only consider line segments of the form $K=[-A,A]$. Considering the previous section, it is not too surprising that the conclusion depends on the $\rank A$.
\begin{prop}
Let $K =[-A,A] \subset \R^{N \times n}$  Suppose that $u \in \WW^{1,2}_\text{loc}(\Omega;\R^N)$ is a local minimiser of $\int_\Omega F_K^2(Dv) \dd x$. 
\begin{enumerate}[(i)]
    \item If $\rank A=1$ then $u \in \WW^{1,\infty}_\mathrm{loc}(\Omega;\R^N)$ .
    \item If $\rank A \geq 2$ then $u \in C^{1,\alpha}_\mathrm{loc}(\Omega;\R^N)$ for some $\alpha>0$.
\end{enumerate}
\end{prop}

\begin{proof}
Firstly we begin with the case $\rank A=1$. After changing coordinates, we may assume $A = \lambda e_1 \otimes e_1$. In the notation of Proposition \ref{prop vectorial regularity decouples} we have $K=K_1 \times \cdots \times K_N$ where $K_1 = [-\lambda,\lambda]e_1$ and $K_\alpha=\left\{ 0 \right\}$ otherwise, so the result follows. The case $\rank A \geq 2$ is more involved. Firstly, note that by rescaling we may assume that $|A|=1$. Define $$\gamma(t)=\begin{cases}-1 & t \leq -1\\
t & -1 <t<1\\
1 & t \geq 1.
    \end{cases}
$$
Then the projection map is given by $\pi(z)=\gamma \left(A \cdot z \right)A$. Therefore, the Euler--Lagrange equation is given by $$-\Delta u=-\Div (\gamma(A \cdot Du)A).$$
By Schauder estimates and the fact $\gamma$ is Lipschitz, if we can show that $A \cdot Du \in C^{0,\alpha}_\text{loc}(\Omega)$, then it follows that $u \in C^{1,\alpha}_\text{loc}(\Omega;\R^N).$ 
To show this, since $u$ has second derivatives we may apply the differential operator $A \cdot \nab$ onto both sides of the Euler--Lagrange equation to get
$$-\Delta(A \cdot Du)=-\Div (A^T A D\gamma(A \cdot Du))$$
which, after applying the chain rule, rearranges to $$-\Div ((I-\gamma'(A \cdot Du)A^TA) D(A \cdot Du))=0.$$
Using $0 \leq \gamma' \leq 1$, the operator is elliptic in the sense that 
$$(I-\gamma'(A \cdot Du)A^TA)\xi \cdot \xi \geq (1-\Vert A \Vert_{\text{op}}^2)|\xi|^2.$$ 
Since $\rank A \geq 2$ and $|A|=1$, it follows that $\Vert A\Vert_{op}<1$. So, $A \cdot Du$ is a solution of an elliptic scalar equation with bounded measurable coefficients. By De Giorgi's theorem this gives $A \cdot Du \in C_\text{loc}^{0,\alpha}(\Omega)$ for some $\alpha>0$.
\end{proof}
\section{AI Declaration}
The author benefitted from the use of ChatGPT-6 Astra in the development of this paper. The main results, including Theorem \ref{Theorem BMO regularity for distance functional}, \ref{Theorem scalar case regularity non convex distance functional}, and \ref{theorem lipschitz regularity for distance functional under supporting ball condition} were developed with a combination of ideas by the author, the author's supervisor, the literature, and ChatGPT-6. A counterexample to Lipschitz regularity, namely Theorem \ref{theorem counterexample lipschitz vectorial first one}, was first produced by ChatGPT after several iterations. The example presented here differs to the original one produced by ChatGPT, and the author's role in this case was to simplify the argument and present it in a readable form. All text and equations were written by hand solely by the author.


\begin{thebibliography}{99}

\bibitem{Asplund1969}
E. Asplund,
\emph{Chebyshev sets in Hilbert space},
Trans. Amer. Math. Soc. \textbf{144} (1969), 235--240.

\bibitem{balashov2024strong}
M.~V.~Balashov and K.~Z.~Biglov,
\emph{The strong convexity supporting condition and the Lipschitz condition for the metric projection},
Math. Notes \textbf{115} (2024), 164--172.

\bibitem{Balashov2012}
M.~V.~Balashov and M.~O.~Golubev,
\emph{About the Lipschitz property of the metric projection in the Hilbert space},
J. Math. Anal. Appl. \textbf{394} (2012), no.~2, 545--551.


\bibitem{BallJames1987}
J. M. Ball and R. D. James,
\textit{Fine phase mixtures as minimizers of energy},
Arch. Ration. Mech. Anal. \textbf{100} (1987),
no. 1, 13--52.

\bibitem{BallKirchheimKristensen2000}
J.~M.~Ball, B.~Kirchheim, and J.~Kristensen,
\emph{Regularity of quasiconvex envelopes},
Calc. Var. Partial Differential Equations \textbf{11} (2000), no.~4, 333--359.


\bibitem{Beck}
L.~Beck,
\emph{Elliptic Regularity Theory: A First Course},
Lecture Notes of the Unione Matematica Italiana, vol.~19,
Springer, Cham, 2016.


\bibitem{Bertsekas2003}
D. P. Bertsekas, A. Nedic, and A. E. Ozdaglar,
\textit{Convex Analysis and Optimization},
Athena Scientific, Belmont, MA, 2003.




\bibitem{Bousquetetal2018}
P.~Bousquet, L.~Brasco, C.~Leone, and A.~Verde,
\emph{On the Lipschitz character of orthotropic $p$-harmonic functions},
Calc. Var. Partial Differential Equations \textbf{57} (2018), no.~3, Paper No.~88.

\bibitem{ChipotEvans1986}
M.~Chipot and L.~C.~Evans,
\emph{Linearisation at infinity and Lipschitz estimates for certain problems in the calculus of variations},
Proc. Roy. Soc. Edinburgh Sect. A \textbf{102} (1986), no.~3--4, 291--303.

\bibitem{DiBenedettoManfredi1993}
E.~DiBenedetto and J.~J.~Manfredi,
\emph{On the higher integrability of the gradient of weak solutions of certain degenerate elliptic systems},
Amer. J. Math. \textbf{115} (1993), no.~5, 1107--1134.

\bibitem{Diening2012}
L.~Diening, P.~Kaplick\'{y}, and S.~Schwarzacher,
\emph{BMO estimates for the $p$-Laplacian},
Nonlinear Anal. \textbf{75} (2012), no.~2, 637--650.

\bibitem{diening2009}
L.~Diening, B.~Stroffolini, and A.~Verde,
\emph{Everywhere regularity of functionals with $\varphi$-growth},
Manuscripta Math. \textbf{129} (2009), no.~4, 449--481.

\bibitem{Diening2011}
L.~Diening, B.~Stroffolini, and A.~Verde,
\emph{Lipschitz regularity for some asymptotically convex problems},
ESAIM Control Optim. Calc. Var. \textbf{17} (2011), no.~1, 178--189.

\bibitem{dolzmann2012}
G.~Dolzmann,
\emph{Regularity of minimizers in nonlinear elasticity---the case of a one-well problem in nonlinear elasticity},
Tech. Mech. \textbf{32} (2012), no.~2--5, 189--194.

\bibitem{DolzmannKristensen2005}
G.~Dolzmann and J.~Kristensen,
\emph{Higher integrability of minimizing Young measures},
Calc. Var. Partial Differential Equations \textbf{22} (2005), no.~3, 283--301.

\bibitem{DolzmannKristensenZhang2013}
G.~Dolzmann, J.~Kristensen, and K.~Zhang,
\emph{BMO and uniform estimates for multi-well problems},
Manuscripta Math. \textbf{140} (2013), no.~1--2, 83--114.

\bibitem{DuzaarMingione2004}
F.~Duzaar and G.~Mingione,
\emph{The $p$-harmonic approximation and the regularity of $p$-harmonic maps},
Calc. Var. Partial Differential Equations \textbf{20} (2004), no.~3, 235--256.

\bibitem{Fletcher2014}
J.~Fletcher and W.~B.~Moors,
\emph{Chebyshev sets},
J. Aust. Math. Soc. \textbf{98} (2015), no.~2, 161--231.

\bibitem{Fonseca1997}
I.~Fonseca and N.~Fusco,
\emph{Regularity results for anisotropic image segmentation models},
Ann. Scuola Norm. Sup. Pisa Cl. Sci. (4) \textbf{24} (1997), no.~3, 463--499.

\bibitem{Fonseca2002}
I.~Fonseca, N.~Fusco, and P.~Marcellini,
\emph{An existence result for a nonconvex variational problem via regularity},
ESAIM Control Optim. Calc. Var. \textbf{7} (2002), 69--95.


\bibitem{Fuchs1994}
M.~Fuchs,
\emph{Regularity for a class of variational integrals motivated by non-linear elasticity},
Asymptot. Anal. \textbf{9} (1994), no.~1, 23--38.

\bibitem{Fuchs1996}
M.~Fuchs,
\emph{Lipschitz regularity for certain problems from relaxation},
Asymptot. Anal. \textbf{12} (1996), no.~2, 145--151.

\bibitem{GiaquintaModica1986}
M.~Giaquinta and G.~Modica,
\emph{Remarks on the regularity of the minimizers of certain degenerate functionals},
Manuscripta Math. \textbf{57} (1986), no.~1, 55--99.

\bibitem{Giusti}
E.~Giusti,
\emph{Direct Methods in the Calculus of Variations},
World Scientific, River Edge, NJ, 2003.

\bibitem{Grafakos2014}
L.~Grafakos,
\emph{Modern Fourier Analysis},
3rd ed., Graduate Texts in Mathematics, vol.~250,
Springer, New York, 2014.

\bibitem{JohnNirenberg1961}
F.~John and L.~Nirenberg,
\emph{On functions of bounded mean oscillation},
Comm. Pure Appl. Math. \textbf{14} (1961), no.~3, 415--426.

\bibitem{Mingione2006}
G.~Mingione,
\emph{Regularity of minima: an invitation to the dark side of the calculus of variations},
Appl. Math. \textbf{51} (2006), no.~4, 355--426.

\bibitem{Rindler}
F.~Rindler,
\emph{Calculus of Variations},
Universitext, Springer, Cham, 2018.

\bibitem{SchevenSchmidt2010}
C.~Scheven and T.~Schmidt,
\emph{Asymptotically regular problems I: higher integrability},
J. Differential Equations \textbf{248} (2010), no.~4, 745--791.

\bibitem{Strong2026}
A.~Strong,
\emph{Sharp regularity results for variational
integrals asymptotic to the $p$-Dirichlet integral},
Preprint (2026)

\bibitem{sverakyan}
V.~\v{S}ver\'{a}k and X.~Yan,
\emph{Non-Lipschitz minimizers of smooth uniformly convex functionals},
Proc. Natl. Acad. Sci. USA \textbf{99} (2002), no.~24, 15269--15276.

\bibitem{Uhlenbeck1977}
K.~Uhlenbeck,
\emph{Regularity for a class of non-linear elliptic systems},
Acta Math. \textbf{138} (1977), 219--240.


\bibitem{Zhang1999}
K. Zhang,
\textit{On various semiconvex relaxations of the squared-distance function},
Proc. Roy. Soc. Edinburgh Sect. A \textbf{129} (1999), 1309--1323.


\bibitem{Zhang1998}
K. Zhang,
\textit{On some quasiconvex functions with linear growth},
J. Convex Anal. \textbf{5} (1998), no. 1, 133--146.


\end{thebibliography}
\end{document}